\documentclass{amsart}
\usepackage{amsmath}
\usepackage{amscd}
\usepackage{amssymb}
\usepackage{amsthm}
\usepackage{amsfonts}
\usepackage{epsfig}
\usepackage{color}

\def\C{\mathbb{C}}
\def\CP{\mathbb{CP}}

\def\Hom{\textup{Hom}}
\def\t{{\mathfrak t}}

\theoremstyle{plain}
\newtheorem{theorem}{Theorem}

\newtheorem{proposition}[theorem]{Proposition}

\newtheorem{example}[theorem]{Example}

\numberwithin{theorem}{section}

\theoremstyle{definition}
\newtheorem{definition}[theorem]{Definition}

\numberwithin{equation}{section}

\title[An introduction to equivariant cohomology and homology]{An introduction to equivariant cohomology and homology, 
following Goresky, Kottwitz, and MacPherson}
\author{Julianna S. Tymoczko}
\address{Department of Mathematics, University of Michigan, 2074 East Hall, Ann Arbor, MI 48109-1109}
\email{tymoczko@umich.edu}

\begin{document}
\begin{abstract}
This paper provides an introduction to equivariant cohomology and homology 
using the approach of Goresky, Kottwitz, and MacPherson.  
When a group $G$ acts suitably on a variety $X$, the equivariant cohomology
of $X$ can be computed using the combinatorial data of 
a skeleton of $G$-orbits on $X$.  We give both a geometric
definition and the traditional definition of equivariant cohomology.  
We include a discussion of the moment map and an algorithm for finding
a set of generators for the equivariant cohomology of $X$.
Many examples and explicit calculations are provided.
\end{abstract}

\maketitle

\section{Introduction}

\begin{figure}[h]
\begin{center}
\scalebox{.8}{\input{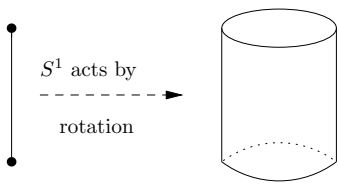}}
\caption{An $S^1$-action} \label{fattening}
\end{center}
\end{figure}

Homology can be thought of as a way of using certain
subspaces of a variety $X$ to describe topological 
invariants of $X$.  Equivariant homology arises in the setting of
a variety $X$ with the action of a group $G$.  Since group
actions typically ``fatten'' a space (see Figure \ref{fattening}), 
information about both the group
and an underlying skeleton in $X$ can be used to recover topological 
information about $X$.

In this paper, we follow the approach of Goresky, Kottwitz, and MacPherson
in \cite{GKM}, who identify a particular skeleton inside $X$ which they
use to turn the problem of computing equivariant cohomology into a
combinatorial puzzle.  The simplest possible skeleton is a point,
but it would be a terribly dull world if the topology of every variety with 
a group action could be reduced to a point.  The next simplest is a graph,
to which the GKM method reduces all suitable varieties.  

Their main result is to show that if the $G$-action on $X$ satisfies certain
conditions, including having a finite number of fixed points and a finite
number of one-dimensional orbits, then the equivariant cohomology of $X$ 
is a 
submodule of the equivariant cohomology of the set of 
fixed points in $X$, which
is naturally a module over a polynomial ring.  
Moreover, this submodule can be identified explicitly using relations
imposed by the one-dimensional orbits.  All of this data can be 
encoded in a graph
whose vertices correspond to fixed points and whose edges correspond to
one-dimensional orbits.

Group actions on varieties had been used to obtain topological data 
long before the paper of Goresky, Kottwitz, and MacPherson.
Bia\l ynicki-Birula used group actions to decompose varieties into a 
collection of attracting cells that could be used to describe cohomology, 
in something like an algebraic analogue of
Morse theory \cite{B-B}.  Working in symplectic geometry, Kirwan and 
many others
had studied symplectic reductions, namely quotients of a variety by its
group action, and had extracted cohomological information using
fixed points and their attracting sets.  \cite[Theorem 5.4]{Ki} 
in fact proves 
 in the symplectic setting that the equivariant cohomology of $X$ is a 
submodule of the equivariant cohomology of the fixed points.  
An early paper of Chang and Skjelbred identifies the image of
the equivariant cohomology of $X$ inside the equivariant cohomology
of the fixed points, for more general topological spaces $X$
\cite[Lemma 2.3]{CS}.  Even earlier contributions 
came in different forms from Borel, Atiyah, 
Hsiang, and Quillen, among others, as described in \cite[Section 1.7]{GKM}.

The combinatorial data associated to group actions on a variety had also 
been studied, especially via the moment map.  The moment map, discussed
at length in Section \ref{moment map}, can be thought of here as a 
combinatorial way to collect data about the group orbits in $X$.
Kostant showed that if a torus $T$ acts on a coadjoint orbit in 
the cotangent space $\t^*$ to $T$, then the image of the moment map is
a convex polytope \cite{Ko}.  This theorem inspired a series of papers
by various authors, including \cite{A}, \cite{H}, and \cite{GuSt} which
extended the result to more general groups acting on more general 
varieties, collectively 
proving that the image of the moment map is some union of
convex polytopes.

What Goresky, Kottwitz, and MacPherson do in their approach to equivariant
cohomology is to put these two schools together, using the combinatorial
data of the orbits of the group action from moment map theory to encode
and interpret topological data of the attracting sets for fixed points
from the study of algebraic group actions.

Many authors have since built on their work.  We mention a sample of 
examples in this introduction and give others at relevant moments in the 
exposition.  GKM theory has been extended to other cohomology theories:
to equivariant intersection homology in \cite{BM} and to equivariant K-theory
in \cite{R(Kn)}.  The class of varieties for which the theory applies has
been expanded, for instance in \cite{GuHo} to 
torus actions without isolated fixed points.  GKM theory is
used to calculate the equivariant cohomology rings of Grassmannians
\cite{KnTao} and regular varieties (which have a unique fixed point with 
respect to a solvable group action) in
\cite{BC}.  The graphs used in GKM theory are studied in their 
own right in \cite{GuZ}.

In order to make this exposition as elementary as possible, we make 
several assumptions throughout this paper.  The first is that we work 
over the complex numbers, so the variety $X$ is a complex algebraic
variety, the group $G$ is a complex linear algebraic group, and all cohomology
has complex coefficients.  The 
second is that the group $G$ is in fact an algebraic torus $T$, so $T$
is a product of $k$ copies of $\C^*$.  When a result is easier to present
using the compact torus $(S^1)^k$, we will use that instead, since the
equivariant cohomology is the same regardless of whether $(\C^*)^k$ or
$(S^1)^k$ is used.
The third assumption is that $X$ is nonsingular and complete.
All of these assumptions can be relaxed, but we make no further mention
of this.  We also note that there are many things in
\cite{GKM} which we do not discuss at all in this paper.  

The rest of this paper provides an introduction to equivariant cohomology
following GKM theory.  Section \ref{defining eq cohom} 
gives several definitions of 
equivariant homology and cohomology.  Section \ref{main theorem} 
then presents the main
theorem due to Goresky, Kottwitz, and MacPherson.  In Section \ref{examples},
 we
give three different examples caculating equivariant cohomology.  We 
include in Section \ref{moment map} 
a discussion of the moment map in general, which
is an important context for the combinatorial graphs used in GKM theory.
Finally, Section \ref{generators} 
describes how to find a minimal set of generators
for the equivariant cohomology,  simplifying the explicit 
calculations that went before.

Many people contributed to the preparation of this
paper through very useful conversations, including
Rebecca Goldin, Allen Knutson, Robert Lazarsfeld, Robert MacPherson, and 
David Nadler.  Tom Braden, Henry Cohn, Allen Knutson, and Robert Lazarsfeld
read drafts of this paper and provided helpful suggestions.  The anonymous
referee's comments were particularly valuable.
The original version of this paper 
was presented the American Mathematical Society summer
conference for young algebraic geometers.  
The author is extremely grateful to the conference organizers
for this opportunity.  

\section{Two definitions of equivariant cohomology and homology}
\label{defining eq cohom}

\subsection{Equivariant homology, topologically defined}\label{top defn}

Let $X$ be an algebraic variety with an algebraic action of
$T = (\C^*)^k$.  To define the ordinary homology of $X$, we
would use the chain group consisting of formal linear combinations of
maps $f:C \longrightarrow X$, where $C$ is a not-too-singular space
whose precise characteristics we do not detail here.  
To incorporate the action of $T$, we modify the definition of the 
chain group as follows:
\[C_i^T(X) = \left\{\begin{array}{c}
\textup{formal linear combinations of $T$-equivariant maps }
f: C \longrightarrow X,\\
\textup{ where $C$ is an $i+k$-dimensional space from a collection } 
\mathcal{C} \\
\textup{of oriented spaces with a free $T$-action and with
a notion of} \\ \textup{ oriented boundaries}\end{array}\right\}.\]
Details about one possible collection $\mathcal{C}$ may be found in 
\cite[Sections 3 and 4]{GKM}.  The notion of oriented boundaries means
that the map $f|_{\partial C}: \partial C \longrightarrow X$ 
is in $C_{i-1}^T(X)$.  Orientations cancel in the usual way, 
so $\partial \partial C=0$.  Consequently, we may take the
homology of the chain complex
\[\cdots C_{i+1}^T(X) \stackrel{\partial_{i+1}}{\longrightarrow} C_i^T(X)
   \stackrel{\partial_{i}}{\longrightarrow} C_{i-1}^T(X) \cdots\]
to get the equivariant homology $H_*^T(X)$.

\begin{figure}[h]
\begin{center}
\input{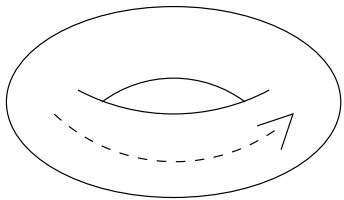}
\caption{The torus with the hula action of $S^1$} \label{hula torus}
\end{center}
\end{figure}

For instance, consider the topological torus with the ``hula-hoop''
action of $S^1$ indicated in Figure \ref{hula torus}.

\begin{figure}[h]
\begin{center}
\input{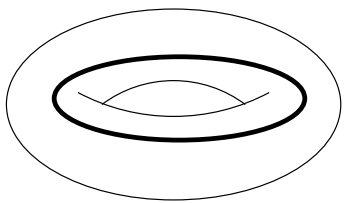}
\caption{The generator for $H_0^{S^1}(X)$} \label{1-chain}
\end{center}
\end{figure}

There is only one free $0+1$-cycle up to homology, shown in Figure
\ref{1-chain}.  Once any given point $p$ on the $0+1$-chain is mapped into
the torus, the $S^1$-equivariance of $f$ determines the rest of the image.  
Any two such cycles $f_1$ and $f_2$ differ by the boundary of a $1+1$-chain,
and so represent the same homology class.
The class of this cycle thus generates $H_0^{S^1}(X)$.

Similarly, there is only one free $1+1$-cycle up to homology.  The image
of this cycle under $f$ must surject onto $X$, 
since any vertical strip in the image 
is wrapped around the torus by $S^1$-equivariance of the map $f$.  Hence,
the class of the cycle given by $X$ itself 
generates $H_1^{S^1}(X)$. 

The reader can observe that the equivariant homology of $X$ with this
$S^1$-action is exactly the same as the ordinary homology of $S^1$.
This is no coincidence.  It is true because the torus is 
$S^1$ (vertically oriented), fattened by a free hula-hoop $S^1$-action.
The more general statement also holds.

\begin{proposition} \label{free action}
If a torus $T$ acts freely on the algebraic variety $X$ and $X/T$
denotes the quotient space of $X$ under this $T$-action, then
$H_*^T(X) = H_*(X/T)$.
\end{proposition}

The proof follows from these definitions.  In fact, $X$ need only be
a principal $T$-bundle for this proposition to hold.

The equivariant cohomology of $X$ may be defined  as the cohomology
of the chain 
complex $\Hom(C_i^T, \C)$.  It is an exercise for the reader that
Proposition \ref{free action} holds for equivariant cohomology, too.

\subsection{Equivariant cohomology, traditionally defined}
The following definition is more common in the literature.  Fix
a $T$-action on $X$.
Proposition \ref{free action} says that equivariant cohomology would
be easy to define if this $T$-action were free.  Consequently, our
strategy is to modify $X$ slightly to obtain a space which does carry a
free $T$-action.

To do this, take $ET$ be a contractible space with a free
$T$-action.  For instance, the group $T=S^1$ acts naturally on each
sphere $S^{2n+1}$ by rotation.  Viewing $S^{2n+1}$ as the unit sphere in
$\C^{n+1}$, the element $e^{i\theta}$ in $S^1$ sends 
$(z_0, z_1, \ldots, z_n)$ to $(e^{i\theta} z_0, e^{i\theta}z_1,
\ldots, e^{i\theta}z_n)$.  This action is free, but unfortunately 
the space $S^{2n+1}$ is not contractible.  To correct for this, observe 
that $S^{2n+1}$ sits inside $S^{2n+3}$ as the equator, and that 
all cycles of dimension at most $2n-2$ are contractible in $S^{2n+3}$.
By taking the union of $S^{2n+1}$ for all $n$, 
we obtain the infinite sphere $S^{\infty}$, which is contractible.  
This space $S^{\infty}$ is $ES^1$.

Define $BT$ to be the quotient $ET/T$; it is called the 
classifying space of $T$.  For example, we may construct the classifying space 
$BS^1$ using the example of $ES^1$ given before.  For each $n$, 
the quotient of $S^{2n+1}$ by the rotation action of $S^1$ is $\CP^{n}$.
Taking the union as before gives a space $\CP^{\infty}$.  Intuitively,
the space $\CP^{\infty}$ is constructed by adjoining cells as if to form
$\CP^n$, but without ever stopping.
%

The diagonal action of $T$ on $X \times ET$ is free, since the action
on $ET$ is free.  Define $X \times_T ET$ to be the quotient 
$(X \times ET) / T$.  This is the total space of a 
fiber bundle with base $BT$ and fiber $X$:
\[\begin{array}{ccc}
X &\longrightarrow & X \times_T ET \\
 & & \downarrow  \\
 & & BT. \\
\end{array}\]
We define the equivariant cohomology of $X$ to be
\[H^*_T(X) = H^*(X \times_T ET).\]

To see what this definition means, we work through the example when $X$ is
a point and $T=S^1$.  
In this case $H^*_T(X) = H^*(\textup{pt} \times_T ET) = 
H^*(ET/T)$.  This is $H^*(BT)$ by definition, and equals 
$H^*(\CP^{\infty})$.
 Recall that
$H^*(\CP^k) = \C[t]/(t^{k+1})$ is generated by the class $[t]$ of 
an embedded $\CP^{k-1}$ together with a normal vector.  The argument does
not generalize precisely, because it uses Poincar\'{e} duality.  However,
the restriction $H^i(\CP^{k+l}) \rightarrow H^i(\CP^k)$ is an isomorphism
as long as $i$ is at most $2k$.
It follows that $H^*(\CP^\infty) = \C[t]$ is the polynomial
ring generated by a single class $[t]$ of complex dimension one.  

Equivariant cohomology satisfies the same canonical properties as ordinary
cohomology.  In particular, the following hold:
\begin{enumerate}
\item functoriality;
\item a ring structure;
\item excision;
\item the Mayer-Vietoris sequence;
\item the K\"{u}nneth formula;
\item the Leray spectral sequence; 
\item for smooth orientable $X$, Poincar\'{e} duality; and
\item existence of Chern classes,
\end{enumerate}
where all subsets and maps are assumed to be equivariant.  

For instance, we can use these properties to show that $H^*_T(X)$ is a
module over $H^*_T(\textup{pt})$.  Applying functoriality to the map
$X \longrightarrow \textup{pt}$ allows us to pull back each class in
$H^*_T(\textup{pt})$ to $H^*_T(X)$.  The ring structure then permits us
to multiply classes coming from $H^*_T(\textup{pt})$ with
classes in $H^*_T(X)$.  

Of course, the same procedure can be followed in 
ordinary cohomology, which brings us to a major difference
between equivariant and ordinary cohomology: $H^*_T(\textup{pt})$
is a very interesting ring, rather than a field.
(One result of this is that the K\"{u}nneth formula and Poincar\'{e}
duality in equivariant cohomology may look slightly unfamiliar.)

To identify the ring $H^*_T(\textup{pt})$, 
suppose the point is acted upon by $T=(S^1)^k$.  
Write $\t$ for the Lie algebra of $T$, namely the tangent to $T$ at 
the identity, and write $\t^*$ for its dual, the cotangent space.  Then
by the same argument as before $H^*_T(\textup{pt}) 
= H^*(\prod_{i=1}^k \CP^\infty)$,
which by the K\"{u}nneth formula is $\otimes_{i=1}^k \C[t_i]$.  We prefer
to write this as the symmetric algebra over the cotangent to $T$,
namely $S(\t^*)$.  (The isomorphism $H^*_T(\textup{pt}) \cong S(\t^*)$
is canonical, though the splitting of $BT = \prod_{i=1}^k \CP^\infty$
need not be.)  This shows that for each finite-dimensional torus
the equivariant cohomology of a point $H^*_T(\textup{pt})$ is $S(\t^*)$.

One way to compute the ordinary cohomology of a fiber bundle
$X \times_T ET$ is to use a Leray spectral sequence.  This is in general
a complicated problem.  Consequently, we restrict to a special class
of varieties $X$ for which this spectral sequence degenerates in the
sense of the next definition.

\begin{definition}
Fix a variety $X$ with an action of $T$.  We say that $X$ 
is equivariantly formal with respect to this $T$-action if
$E^2=E^\infty$ in the spectral sequence associated to the fibration
$X \longrightarrow X \times_T ET \longrightarrow BT$.
\end{definition}
 
In algebraic terms, equivariant formality implies that $H^*_T(X)$ is a free 
module over $H^*_T(\textup{pt})$.  We will use this crucial fact 
repeatedly in what follows.
The next proposition follows directly from the definition and is the
reason we consider equivariant formality.

\begin{proposition} \label{point of eq formal}
If $X$ is equivariantly formal with respect to a given $T$-action then
there is a natural $H^*_T(\textup{pt})$-module 
isomorphism $H^*_T(X) = H^*(X) \otimes H^*_T(\textup{pt})$.
\end{proposition}

When $X$ is equivariantly formal with respect to $T$, 
the ordinary cohomology of $X$ can be reconstructed from its
equivariant cohomology.  To do this, let $M$ be the augmentation ideal
in $S(\t^*)$.  (Recall that $M$ is the maximal homogeneous
ideal in $S(\t^*)$, so if $\t^*$ is spanned by basis elements
$t_1$, $t_2$, $\ldots$, $t_k$ then $M = 
\langle t_1, \ldots, t_k \rangle$.)  The ordinary cohomology
can be reconstructed as the quotient
\[H^*(X) = \frac{H^*_T(X)}{M \cdot H^*_T(X)}\]
which in effect simply sets each $t_i=0$.  This
is an equality of rings and not just of modules \cite[1.2.4]{GKM}.  
Whenever we compute an example of $H^*_T(X)$, 
we will also give generators
for $H^*_T(X)$ as a module over $H^*_T(\textup{pt})$ that coincide with
module generators for $H^*(X)$.

In practice, many varieties of interest are equivariantly formal.  For
instance, \cite{GKM} note that all of the following are equivariantly 
formal:
\begin{enumerate}
\item a smooth complex projective 
   algebraic variety (with respect to any linear algebraic 
   $T$-action);
\item a variety whose ordinary cohomology vanishes in odd degree (with
  respect to any $T$-action); 
\item a variety with a $T$-invariant CW-decomposition; and
\item a compact symplectic manifold with a Hamiltonian $T$-action, 
   where $T$ is a compact torus.
\end{enumerate}

These are listed in \cite[Section 1.2 and Theorem 14.1]{GKM}.

\begin{figure}[h]
\begin{center}
\scalebox{.4}{\input{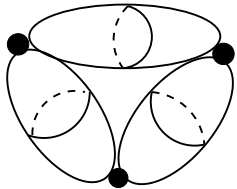}}
\caption{A complete variety which is not equivariantly formal} \label{non equivariantly formal}
\end{center}
\end{figure}

The quintessential non-equivariantly formal space is shown in Figure
\ref{non equivariantly formal}; the
points of intersection are poles of the spheres.  
The relevant action of $\C^*$ fixes the points of intersection and 
rotates each copy of $\CP^1$ so that when acted on by $t$ in $\C^*$,
as $t \rightarrow 0$ the points 
of each $\CP^1$ flow from the north pole to the south pole, and so that
the north pole of each sphere is glued to the south pole of the next.

\section{The main theorem} \label{main theorem}

We now come to the main construction in this paper.
Suppose $X$ has a $T$-action with respect to which three conditions hold:
\begin{enumerate}
\item $X$ is equivariantly formal;
\item $X$ has finitely many fixed points; and
\item $X$ has finitely many one-dimensional orbits.
\end{enumerate}
In this situation, Goresky, Kottwitz, and MacPherson show that
the combinatorial data encoded in the graph of fixed points
and one-dimensional orbits of $T$ in $X$
implies a particular algebraic characterization of $H^*_T(X)$.  For the
rest of this paper, we will use the algebraic torus $(\C^*)^k$ instead of 
its deformation retract $(S^1)^k$.

Given a space $X$ with a $T$-action, the graph
of fixed points 
and one-dimensional orbits can often be computed directly.  For instance, let
$X$ be $\CP^2$, written with homogeneous coordinates.
Define an action of $(\C^*)^2$ on $X$ by $(t_1, t_2) \cdot
[x_0, x_1, x_2] = [x_0, t_1x_1, t_2x_2]$.  The point $[x_0,x_1,x_2]$
lies in a two-dimensional orbit unless it is zero in one or two
coordinates.  Thus, there are exactly three fixed
points under this action, namely $[1,0,0]$, $[0,1,0]$, and $[0,0,1]$,
and exactly three one-dimensional orbits, namely $[*,*,0]$, $[0,*,*]$, 
and $[*,0,*]$, where $*$ denotes any nonzero element in $\C$.
The closure of each one-dimensional orbit is isomorphic to $\CP^1$ and has 
exactly two fixed points in its closure, one at the north pole and one at
the south pole.
Note that 
the fixed points and
one-dimensional orbits form a subvariety isomorphic 
to that
of Figure \ref{non equivariantly formal}, although the $T$-action is
different.  

In fact, the basic properties of this subvariety hold for the subvariety
of fixed points and one-dimensional orbits of any variety $X$ satisfying the
same hypotheses.

\begin{proposition}
Let $X$ be a complete variety with a fixed action of $T$, with respect to which
$X$ is equivariantly formal, and which has 
finitely many fixed points and finitely 
many one-dimensional orbits.  The following all hold:
\begin{enumerate}
\item The closure of each one-dimensional orbit of $T$ in $X$ is
isomorphic to $\CP^1$.  
\item There are exactly two fixed points in the closure
of each one-dimensional orbit, one at the north pole and the other at the south pole.
\item The torus $T$ acts on each one-dimensional orbit by rotation.
\item For each one-dimensional orbit $O$, there is a subtorus $T' \subseteq T$ of 
codimension one which fixes $O$ pointwise.
\end{enumerate}
\end{proposition}

These facts can all be found in \cite[Section 7.1]{GKM}.

This proposition implies that the fixed points and one-dimensional orbits 
of such a
variety can be represented by a graph, each of whose vertices corresponds
to a fixed point and each of whose edges corresponds to a one-dimensional 
orbit whose
poles are the endpoints of the edge.  This graph is called the moment graph 
of $X$, discussed further in Section \ref{moment map}.  Note that
the proposition further states that each edge is associated to a unique
subtorus of codimension one that fixes the corresponding one-dimensional 
orbit pointwise.  We typically label each edge with the annihilator in
the cotangent space of the tangent to this stabilizer.  
With our assumptions on $X$, the edges adjacent to a 
given vertex will be distinct.

There is a natural inclusion from the set of $T$-fixed points in $X$ to
all of $X$.  This map induces a map on equivariant cohomology 
$H^*_T(X) \longrightarrow \bigoplus_{p \textup{ is $T$-fixed}}
H^*_T(p)$.  From what we said earlier, it follows that $H^*_T(X)$ maps to
$\bigoplus_{p \textup{ is $T$-fixed}}
H^*_T(p) = \bigoplus_{p \textup{ is $T$-fixed}} S(\t^*)$.
Given our hypotheses on $X$ and the $T$-action, this map is actually
an injection.

\begin{theorem} \label{injection theorem}
Suppose $X$ carries a $T$-action with respect to which $X$ is
equivariantly formal, and which has finitely many fixed points.
Then the map
\[H^*_T(X) \longrightarrow \bigoplus_{p \textup{ is $T$-fixed}} S(\t^*)\]
induced by inclusion from the $T$-fixed points into $X$ is an injection.
\end{theorem}

This is proven in \cite{GKM}, but was also proven previously in 
\cite[Theorem 5.4]{Ki}.
There are earlier results that 
anticipate this, including those of Borel, Atiyah, Hsiang, and Quillen, 
listed along with a more thorough historical
bibliography in \cite[Section 1.7]{GKM}.
A particularly pretty proof in the case
of symplectic manifolds is given in \cite{ToW}.

The reader is encouraged to treat the next example as an exercise, though
a sketch of the proof  is included.

\begin{example}
The variety $X$ in Figure \ref{non equivariantly formal} is not
equivariantly formal.
\end{example}

\begin{proof}
Denote the $T$-fixed points of $X$ by $X^T$.  By excision, the group
$H^1_T(X,X^T)$ is three-dimensional.  The long exact sequence for the
pair $(X,X^T)$ shows that $H^1_T(X)$ does not vanish.  If $X$ were 
equivariantly formal, then $H^*_T(X) \longrightarrow H^*_T(X^T)$ would
be an injection, but $H^*_T(X^T)$ does not exist in odd degrees.
\end{proof}

We wish to identify the image of $H^*_T(X)$ inside the ring
$\bigoplus_{p \textup{ is $T$-fixed}} S(\t^*)$.  Suppose $O$
is a one-dimensional orbit with north pole $N$ and south pole $S$, 
and let $f_N$ and 
$f_S$ be elements of $S(\t^*)$ associated to the north and south pole
respectively.  There is a subtorus $T' \subseteq T$ of codimension one
which fixes $O$.  We wish to see that in order to be in the image of 
$H^*_T(X)$, the polynomial functions associated to the poles
of $O$ must agree on the tangent space $\t'$ of $T'$, namely
the restrictions $f_N|_{\t'} = f_S|_{\t'}$.  Equivalently, we need
to show that the 
difference $f_N - f_S$ lies in the ideal generated by the 
annihilator of $\t'$ inside the cotangent $\t^*$.  The 
following commutative diagram, all of whose maps are equivariant, proves
that this claim is true.  (The 
identification of $H^*_T(O)$ is implied by the construction in Section
\ref{top defn}.)
\[\begin{array}{ccccccc}
{} && {H^*_T(\overline{O})} & \longrightarrow& {H^*_T(O \cup N)} & \cong & 
   {H^*_T(N)} \\
 && \downarrow & & \downarrow && \\
{H^*_T(S)} &\cong & {H^*_T(O \cup S)} &\longrightarrow & {H^*_T(O)} &
\cong & {S((\t')^*)}
\end{array} \]

In fact, the condition that $f_N|_{\t'} = f_S|_{\t'}$ 
is sufficient to identify the image of the map
$H^*_T(X) \longrightarrow \bigoplus_{p \textup{ is $T$-fixed}} 
S(\t^*)$.  This is the main content of the following 
theorem of Goresky, Kottwitz, and MacPherson.

\begin{theorem} \label{GKM theorem}
Let $X$ be an algebraic variety with a $T$-action with respect to which
$X$ is equivariantly formal, and which has 
finitely many fixed points and finitely
many one-dimensional orbits.  Denote the one-dimensional orbits $O_1$, $\ldots$, $O_m$.  For each
$i$, denote the poles of $O_i$ by $N_i$ and $S_i$ and denote the 
stabilizer of $O_i$ in $T$ by $T_i$.
The equivariant cohomology ring of $X$ is given by 
\[H^*_T(X) \cong \left\{(f_{p_1}, \ldots, f_{p_m}) 
\in \bigoplus_{\textup{fixed pts}} S(\t^*):  f_{N_i}|_{\t_i} = 
f_{S_i}|_{\t_i} \textup{ for each }i=1,\ldots,m \right\}.\]
\end{theorem}

We stress that this is a ring isomorphism, unlike the isomorphism in
Proposition \ref{point of eq formal}.  That the relations $f_{N_i}|_{\t_i} = 
f_{S_i}|_{\t_i}$ suffice to determine the image of $H^*_T(X)$ was proven
earlier in \cite[Lemma 2.3]{CS}, though at a level of generality which
prevented the concrete description which follows.  
We also remark that \cite{GuHo}
contains a very nice sketch of the proof of this theorem in the symplectic
setting, assuming the injection of Theorem \ref{injection theorem}.

\begin{figure}[h]
\begin{picture}(300,15)(0,-5)
\put(100,0){\circle*{5}}
\put(200,0){\circle*{5}}
\put(100,0){\line(1,0){100}}
\put(150,5){$0$}
\put(70,0){$[1,0]$}
\put(210,0){$[0,1]$}
\end{picture}
\caption{The moment map graph for $\CP^1$} \label{moment for CP1}
\end{figure}
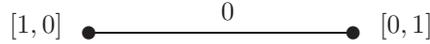

For instance, let $X$ be $\CP^1$ with an action of $T=\C^*$ given by
$t \cdot [x_0,x_1] = [x_0,tx_1]$.  There are exactly two fixed points
of this action, namely $[1,0]$ and $[0,1]$, and there is exactly one
one-dimensional orbit given by $[*,*]$.  The graph associated to these fixed
points and one-dimensional orbits is given in Figure \ref{moment for CP1}.  
Since
the subtorus which fixes the one-dimensional orbit is simply the 
identity, its
tangent space is $\{0\}$, with which we have labelled the edge in
Figure \ref{moment for CP1}.  We use the polynomial description of the
equivariant cohomology of each fixed point.
The polynomials associated to each
fixed point must agree on this tangent space, so $p_N(0)=p_S(0)$.  By
Theorem \ref{GKM theorem} (and used in its proof), 
the equivariant cohomology of $\CP^1$ is
\[H^*_{\C^*}(\CP^1) = \{(p_N,p_S) \in S(\t^*) \oplus S(\t^*):
p_N(0)=p_S(0)\}.\]

To verify this answer, we first rewrite the condition on the polynomials
associated to the fixed points.  If $p_N$ is chosen
freely from $S(\t^*)$ then $p_S = p_N + tp$ for a polynomial $p$ freely
chosen in $S(\t^*)$.  In other words, the equivariant cohomology of
$\CP^1$ is generated additively by the classes $(p_N,p_N)$ and $(0,tp)$,
where $p_N$ and $p$ run over an additive basis of $S(\t^*)$.
Using Proposition \ref{point of eq formal}, we see that 
$H^*_{\C^*}(\CP^1) = H^*(\CP^1) \otimes H^*_{\C^*}(\textup{pt})$.  In
this case, we know that the ordinary cohomology $H^*(\CP^1)$ is generated
by a class $[1]$ in degree zero and $\alpha$ in degree $2$.  Our calculations
identified one family $(p_N,p_N)$ corresponding to 
$[1] \otimes S(\t^*)$
and a second family $(0,tp)$ corresponding to $\alpha \otimes S(\t^*)$, so
our computation of $H^*_{\C^*}(\CP^1)$ is confirmed.

\section{Examples} \label{examples}

This section contains more examples computing equivariant cohomology.
The calculations here are ad hoc and follow directly from the definitions
and theorems in the previous sections.  Section \ref{generators} gives
a more methodical way to construct generators for $H^*_T(X)$, systematizing
some of these calculations.  The reader is encouraged to treat these 
examples as exercises.

\begin{example} \label{cp2 example}
The equivariant cohomology of $\CP^2$ with an action of
$(\C^*)^2$.
\end{example}

\begin{figure}[h]
\begin{picture}(300,40)(0,0)
\put(100,0){\circle*{5}}
\put(200,0){\circle*{5}}
\put(150,20){\line(-5,-2){50}}
\put(150,20){\line(5,-2){50}}
\put(150,20){\circle*{5}}
\put(100,0){\line(1,0){100}}
\put(135,35){$[1,0,0]$}
\put(145,25){$p_1$}
\put(68,10){$[0,1,0]$}
\put(85,0){$p_2$}
\put(210,10){$[0,0,1]$}
\put(210,0){$p_3$}
\put(115,15){$t_1$}
\put(180,15){$t_2$}
\put(135,-10){$t_1-t_2$}
\end{picture}
\caption{The fixed points and one-dimensional orbits for $\CP^2$} \label{cp2}
\end{figure}
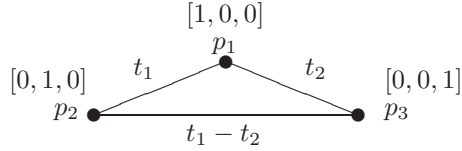

If $[x_0, x_1, x_2]$ is a point of $\CP^2$ and $(t_1,t_2)$ is an 
element of the torus $(\C^*)^2$, then the torus action is given by
$(t_1,t_2)\cdot [x_0, x_1, x_2] = [x_0, t_1 x_1, t_2x_2]$.  As 
described in Section \ref{main theorem}, there are three fixed points
and three one-dimensional orbits under this action.  The one-dimensional orbit whose closure is
$[*,0,*]$ is stabilized by the subtorus $\C^* \times 1$ while
the one-dimensional orbit whose closure is
$[*,*,0]$ is stabilized by the subtorus $1 \times \C^*$.  Finally, the
one-dimensional orbit $[0,*,*]$ is stabilized by the one-dimensional subtorus given
by $(t,t)$ for $t \in \C^*$.  Figure \ref{cp2} shows the graph of
fixed points 
and one-dimensional orbits for this function, labelling each edge 
according to
the annihilator in the cotangent space of the stabilizer of the orbit.
(Note that we use $t_i$ to denote the standard basis of $\t^*$, which
is dual to that of $\t$.)

We now analyze the equivariant cohomology ring using Theorem 
\ref{GKM theorem}.  Suppose $p_1$, $p_2$, and $p_3$ are elements of
$S(\t^*)$ associated to the vertices as 
labelled in Figure \ref{cp2}.  Then the theorem tells us that
\[\begin{array}{rcl}
p_2-p_1 & \in & \langle t_1 \rangle, \\
p_3-p_1 & \in & \langle t_2 \rangle, \textup{ and } \\
p_3-p_2 & \in & \langle t_2 - t_1 \rangle. \end{array}\]
These relations show that we may choose arbitrary elements $p_1$ and $q_1$ 
of $S(\t^*)$ and set $p_2 = p_1 + t_1q_1$.  Once we have done that, 
we  choose an element $q_2$ 
of $S(\t^*)$ and determine the conditions imposed on $q_2$ by setting 
$p_3 = p_1 + t_2q_2$.  Subtracting shows that
$p_3 - p_2 = t_2 q_2 - t_1 q_1$.  Since this is in the
ideal $\langle t_2 - t_1 \rangle$, we may rewrite $q_2 = q_1 + (t_2-t_1)q$.
This equation shows that the most general choice of $q_2$ in 
$S(\t^*)$ corresponds to a free choice of $q$
in $S(\t^*)$.  In other words, the equivariant
cohomology ring of $\CP^2$ consists of the triples
\[H_{(\C^*)^2}^*(\CP^2) = \{(p_1, p_1 + t_1 q_1, p_1 + t_2q_1 + t_2(t_2-t_1)q): 
   p_1, q_1, q \in S(\t^*)\}.\]
Note that the equivariant cohomology $H_{(\C^*)^2}^*(\CP^2)$ 
is generated as a module over $S(\t^*)$ by the elements $(1, 1, 1)$,
$(0,t_1,t_2)$, and $(0,0,t_2(t_2-t_1))$ of degrees $0$, $1$, and $2$
respectively.  This confirms  our calculation, since  
$H_{(\C^*)^2}^*(\CP^2) = H^*(\CP^2) \otimes H_{(\C^*)^2}^*(\textup{pt})$ by
Proposition \ref{point of eq formal}, which simplifies to 
$H^*(\CP^2) \otimes S(\t^*)$.

\begin{example} \label{flag variety example}
The equivariant cohomology of the full flag variety in $\C^3$ 
under the diagonal action of $(\C^*)^2$.
\end{example}

In this example we consider the full flag variety in $\C^3$, namely the
collection of nested vector spaces $V_1 \subsetneq V_2 \subsetneq V_3 = \C^3$
with each $V_i$ an $i$-dimensional subspace.  To represent the flag
$V_1 \subsetneq V_2 \subsetneq V_3$, we choose vectors $g_1$, $g_2$, and
$g_3$ so that $V_i$ is the linear span of $g_1$ through $g_i$.  We then
represent the flag by the matrix whose $i^{th}$ column vector is $g_i$.  
(This representation is not unique.)

Let $g$ be a matrix whose entries are $g_{i,j}$ and let $(t_1, t_2)$
be an element of $\C^* \times \C^*$.  The element $(t_1,t_2)$ acts
on $g$ by conjugation, so that 
$(t_1, t_2) \cdot g$ is given by
\[\left( \begin{array}{ccc} t_1 & 0 & 0 \\ 0 & 1 & 0 \\ 0 & 0 & t_2 
\end{array} \right) g \left( \begin{array}{ccc} \frac{1}{t_1} & 0 & 0 \\ 
0 & 1 & 0 \\ 0 & 0 & \frac{1}{t_2} \end{array} \right)\]
though multiplication on the right by a diagonal matrix fixes each flag.

Suppose $g$ is a matrix representing a flag.  If the flag given by $g$ 
is fixed by this
torus action, then the first column of $g$ 
can only have one nonzero entry $g_{w_1,1}$.
We may assume that the entries in the columns to the right of $g_{w_1,1}$ are
zero by choosing $g_2$ and $g_3$ appropriately.
By the same argument, the second column of $g$ can only have one nonzero
entry $g_{w_2,2}$ 
and we may assume that the entry $g_{w_2,3}$ is zero without
altering the flag.  This leaves a unique entry $g_{w_3,3}$ nonzero.
In sum, the flags fixed by this torus action are exactly the flags generated
by permutations of the basis vectors.

A similar argument shows that the one-dimensional orbits are parametrized by 
the matrices which are a permutation matrix plus a single nonzero entry $a$.
Moreover, the entry $a$ must lie to the left and above the nonzero
entries of the permutation matrix.  A set of matrices parametrizing each
 one-dimensional orbit is given in Figure
\ref{1-orbits of flag}.

\begin{figure}[h]
\[\begin{array}{cccc}
\left( \begin{array}{lll} 0 & 0 & 1 \\ a & 1 & 0 \\ 1 & 0 & 0 \end{array} \right) &
\multicolumn{2}{c}{\left( \begin{array}{lll} a & 0 & 1 \\ 0 & 1 & 0 \\ 1 & 0 & 0 \end{array}
 \right)} &
\left( \begin{array}{lll} 0 & a & 1 \\ 0 & 1 & 0 \\ 1 & 0 & 0 \end{array} \right) \\ \\
\left( \begin{array}{lll} 0 & a & 1 \\ 1 & 0 & 0 \\ 0 & 1 & 0 \end{array} \right) &
\left( \begin{array}{lll} a & 0 & 1 \\ 1 & 0 & 0 \\ 0 & 1 & 0 \end{array} \right) &
\left( \begin{array}{lll} 0 & 1 & 0 \\ a & 0 & 1 \\ 1 & 0 & 0 \end{array} \right) &
\left( \begin{array}{lll} a & 1 & 0 \\ 0 & 0 & 1 \\ 1 & 0 & 0 \end{array} \right) \\ \\
\multicolumn{2}{c}{\left( 
\begin{array}{lll} a & 1 & 0 \\ 1 & 0 & 0 \\ 0 & 0 & 1 \end{array} \right)} &
\multicolumn{2}{c}{\left( \begin{array}{lll} 1 & 0 & 0 \\ 0 & a & 1 \\ 0 & 1 & 0
\end{array} \right)}
\end{array}\]
\caption{Matrices parametrizing one-dimensional orbits of the flag variety in $\C^3$, for $a$ in $\C$} \label{1-orbits of flag}
\end{figure}

Finally, we identify the closures of these one-dimensional orbits.  When $a$ is 
zero, the permutation is clear.  When $a$ goes to $\infty$, the one-dimensional orbit
approaches a transposition of this permutation.  We give an example; all
the rest follow the same pattern.
\[\begin{array}{lcl}
\left( \begin{array}{ccc} a & 0 & 1 \\ 0 & 1 & 0 \\ 1 & 0 & 0 
\end{array}\right) &= &\textup{ Flag }
\langle ae_1 + e_3 \rangle \subsetneq
  \langle ae_1 + e_3, e_2 \rangle \subsetneq
  \langle e_1, e_2, e_3 \rangle \\
&= &\textup{ Flag }
\langle e_1 + a^{-1}e_3 \rangle \subsetneq
  \langle e_1 + a^{-1}e_3, e_2 \rangle \subsetneq
  \langle e_1, e_2, e_3 \rangle \vspace{.15in}\\
& \stackrel{a \to
  \infty}{\longrightarrow} &\textup{ Flag }
\langle e_1 \rangle \subsetneq
  \langle e_1, e_2 \rangle \subsetneq
  \langle e_1, e_2, e_3 \rangle  =  
\left( \begin{array}{ccc} 1 & 0 & 0 \\ 0 & 1 & 0 \\ 0 & 0 & 1 
\end{array}\right) \end{array}\]
(In fact, the two fixed points in the closure of each one-dimensional
orbit are always related by multiplication by a transposition.  We do not
give a systematic description of this in order to maintain an elementary
exposition while avoiding cumbersome notation.)

The stabilizer of each one-dimensional orbit can be determined from Figure
\ref{1-orbits of flag} using the multiplication rule.
This is all of the information needed to build the graph of fixed points and
one-dimensional orbits for
this torus action on the flag variety, shown in Figure \ref{moment map
for flags}.  The edges have been arranged so that they correspond to the
matrices in Figure \ref{1-orbits of flag}, ordered from left to right 
and top to bottom according to the
endpoint which is highest on the graph.  Each edge is labelled by the
linear functional in $\t^*$ that annihilates that edge's stabilizer in $T$.
Note that if two edges in this graph are parallel then they
share the same label.

\begin{figure}[h]
\begin{picture}(350,120)(0,-60)
\put(175,45){\circle*{5}}
\put(175,-45){\circle*{5}}
\put(125,20){\circle*{5}}
\put(225,20){\circle*{5}}
\put(125,-20){\circle*{5}}
\put(225,-20){\circle*{5}}
\put(175,45){\line(-2,-1){50}}
\put(175,45){\line(2,-1){50}}
\put(175,45){\line(0,-1){90}}
\put(125,20){\line(0,-1){40}}
\put(125,20){\line(5,-2){100}}
\put(225,20){\line(0,-1){40}}
\put(225,20){\line(-5,-2){100}}
\put(125,-20){\line(2,-1){50}}
\put(225,-20){\line(-2,-1){50}}
\put(150,-60){\tiny $\left( \begin{array}{ccc} 1 & 0 & 0 \\
      0 & 1 & 0 \\ 0 & 0 & 1 \end{array} \right)$}
\put(150,55){\tiny $\left( \begin{array}{ccc} 0 & 0 & 1 \\
      0 & 1 & 0 \\ 1 & 0 & 0 \end{array} \right)$}
\put(70,30){\tiny $\left( \begin{array}{ccc} 0 & 0 & 1 \\
      1 & 0 & 0 \\ 0 & 1 & 0 \end{array} \right)$}
\put(225,-35){\tiny $\left( \begin{array}{ccc} 1 & 0 & 0 \\
      0 & 0 & 1 \\ 0 & 1 & 0 \end{array} \right)$}
\put(225,30){\tiny $\left( \begin{array}{ccc} 0 & 1 & 0 \\
      0 & 0 & 1 \\ 1 & 0 & 0 \end{array} \right)$}
\put(70,-35){\tiny $\left( \begin{array}{ccc} 0 & 1 & 0 \\
      1 & 0 & 0 \\ 0 & 0 & 1 \end{array} \right)$}
\put(145,-40){$t_1$}
\put(200,-40){$t_2$}
\put(95,0){$t_1-t_2$}
\put(230,0){$t_1-t_2$}
\put(145,37){$t_2$}
\put(200,37){$t_1$}
\put(145,-7){$t_2$}
\put(200,-7){$t_1$}
\put(177,20){$t_1-t_2$}
\end{picture}
\caption{The graph of fixed points and one-dimensional orbits for the flag variety}
\label{moment map for flags}
\end{figure}
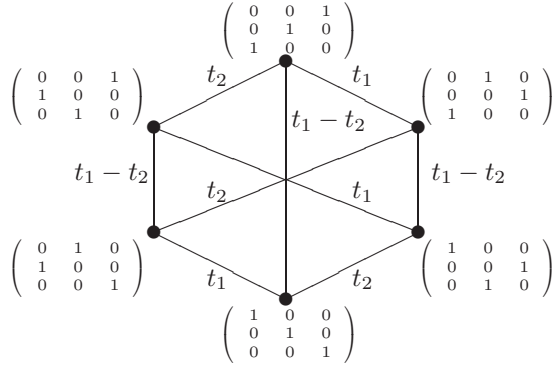

\begin{figure}[h]
\begin{picture}(350,135)(0,-60)
\put(175,45){\circle*{5}}
\put(175,-45){\circle*{5}}
\put(125,20){\circle*{5}}
\put(225,20){\circle*{5}}
\put(125,-20){\circle*{5}}
\put(225,-20){\circle*{5}}
\put(175,45){\line(-2,-1){50}}
\put(175,45){\line(2,-1){50}}
\put(175,45){\line(0,-1){90}}
\put(125,20){\line(0,-1){40}}
\put(125,20){\line(5,-2){100}}
\put(225,20){\line(0,-1){40}}
\put(225,20){\line(-5,-2){100}}
\put(125,-20){\line(2,-1){50}}
\put(225,-20){\line(-2,-1){50}}
\put(145,-40){$t_1$}
\put(200,-40){$t_2$}
\put(95,0){$t_1-t_2$}
\put(230,0){$t_1-t_2$}
\put(145,37){$t_2$}
\put(200,37){$t_1$}
\put(145,-7){$t_2$}
\put(200,-7){$t_1$}
\put(177,20){$t_1-t_2$}
\put(170,-55){$p_1$}
\put(80,-25){$p_1 +t_1 p_2$}
\put(230,-25){$p_1 + t_2 p_3$}
\put(30, 25){$\begin{array}{l} p_1 + t_1 p_2 - (t_1-t_2) p_3 \\
    + (t_1-t_2) t_1 p_4 \end{array}$}
\put(220, 25){$\begin{array}{r} p_1 + (t_1-t_2) p_2 + t_2 p_3 \\
    + (t_1-t_2) t_2 p_5\end{array}$}
\put(110, 65){$\begin{array}{c} p_1 + (t_1-t_2) p_2 - (t_1-t_2) p_3 \\
    + (t_1-t_2) t_1 p_4+(t_1-t_2) t_2 p_5 \\
   + t_1 t_2 (t_1-t_2)p_6 \end{array}$}
\end{picture}
\caption{Functions in $S(\t^*)$ associated to each vertex}
\label{eqcohom of flags}
\end{figure}
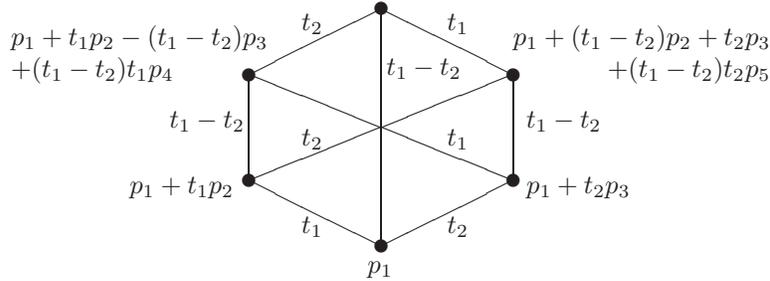

Now we use Theorem \ref{GKM theorem} to compute the equivariant 
cohomology of the flag variety.  We represent this calculation 
in Figure \ref{eqcohom of flags} using
the previous graph, labelling each vertex this time with an element of
$S(\t^*)$.  
If $p_1$ is the function in $S(\t^*)$ 
associated to the lowest vertex, then the two vertices immediately
above the lowest vertex
can be labelled as indicated in Figure \ref{eqcohom of flags}, 
for any $p_2$ and $p_3$ in $S(\t^*)$.
The two vertices above those have conditions imposed by the cycles in the 
graph.  For instance, the vertex on the left can be written as
$p_1 + t_1 p_2 + (t_1-t_2) q$ by moving up the edge labelled $t_1-t_2$
from level one.  Since this function can also be obtained by moving over
the edge labelled $t_1$ from level one, we conclude that 
\[(p_1 + t_1 p_2 + (t_1-t_2) q) - (p_1 + t_2 p_3) \in 
\langle t_1 \rangle.\]
In other words, the difference $t_2 (q + p_3)$ is in $\langle t_1 \rangle$
and so $q = t_1 p_4 - p_3$ for $p_4$ chosen freely in $S(\t^*)$.  
We have labelled this vertex accordingly in Figure \ref{eqcohom of flags}.  
The calculation for the vertex on the right follows a symmetric
argument.  The
topmost vertex must be labelled with a function specified by three edge
relations, given in Figure \ref{eqcohom of flags}.  

\begin{figure}[h]
\begin{picture}(350,50)(0,-25)
\put(130,25){$(0,0,0,0,0,t_1t_2(t_1-t_2))$}
\put(15,8){$(0,0,0,t_1(t_1-t_2),0,t_1(t_1-t_2))$}  
\put(195,8){$(0,0,0,0,t_2(t_1-t_2),t_2(t_1-t_2))$}
\put(195,-12){$(0,0,t_2,-t_1+t_2, t_2, -t_1+t_2)$}
\put(40,-12){$(0, t_1, 0, t_1, t_1-t_2, t_1-t_2)$}
\put(145,-33){$(1,1,1,1,1,1)$}
\put(175,20){\circle*{5}}
\put(165,10){\circle*{5}}
\put(185,10){\circle*{5}}
\put(165,-10){\circle*{5}}
\put(185,-10){\circle*{5}}
\put(175,-20){\circle*{5}}
\put(175,20){\line(-1,-1){10}}
\put(175,20){\line(1,-1){10}}
\put(175,20){\line(0,-1){40}}
\put(165,10){\line(0,-1){20}}
\put(165,10){\line(1,-1){20}}
\put(185,10){\line(0,-1){20}}
\put(185,10){\line(-1,-1){20}}
\put(175,-20){\line(-1,1){10}}
\put(175,-20){\line(1,1){10}}
\end{picture}
\caption{Generators for the equivariant cohomology of the flag variety}
\label{flag generators}
\end{figure}
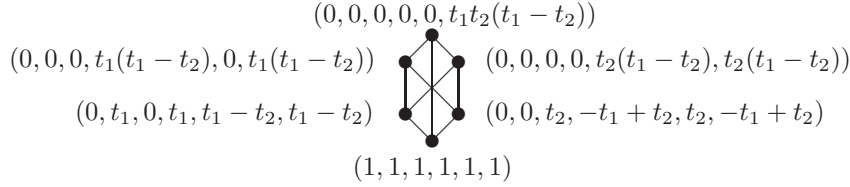

In other words, the equivariant cohomology of the flag variety is generated
as an $S(\t^*)$-module by six elements.  We list these generators in
Figure \ref{flag generators}, associating each generator to
the vertex where the function first appears.
Figure \ref{generator examples} 
displays several of these generators in a style used
frequently by those who work in this area.
Note that as a module over $S(\t^*)$ the equivariant cohomology has 
one generator of dimension zero, two of dimension one,
two of dimension two, and one of dimension three.  These are the even
dimensional Betti numbers of the flag variety, as Proposition \ref{point 
of eq formal}
assures us they should be.

\begin{example} \label{G24 example}
The hypersurface $x_1y_1 + x_2y_2+x_3y_3=0$ in $\CP^5$.
\end{example}

We remark that this hypersurface is in fact the
image of the Grassmannian $Gr(2,4)$ under the Plucker embedding.

The torus $(\C^*)^3$ acts coordinatewise 
on this hypersurface as follows.  Given an element 
 $(x_1, x_2, x_3,y_1,y_2,y_3)$
of the hypersurface, then
\[\begin{array}{l}
(t_1,1,1) \cdot (x_1, x_2, x_3,y_1,y_2,y_3) = 
   (t_1x_1, t_1x_2, t_1x_3,y_1,y_2,y_3)\\
(1,t_2,1) \cdot (x_1, x_2, x_3,y_1,y_2,y_3) = 
   (t_2x_1, t_2x_2, x_3,y_1,y_2,t_2y_3) \\
(1,1,t_3) \cdot (x_1, x_2, x_3,y_1,y_2,y_3) = 
   (t_3x_1, x_2, t_3x_3,y_1,t_3y_2,y_3). \end{array}\]

We now verify that each fixed point under this action has exactly 
one nonzero coordinate.  In fact, at most one of the $x_i$ is nonzero,
since otherwise the action of the subtorus $1 \times \C^* \times \C^*$ is at
least one-dimensional.  Likewise, at most one of the $y_i$ is nonzero,
using the action of the same subtorus.  Finally, either one of the $x_i$ or
one of the $y_j$ is nonzero, but not both, since otherwise 
the action of $\C^* \times
1 \times 1$ would form a one-dimensional orbit.
This means that the six coordinate entries are the fixed points of this
torus action on the hypersurface.

The one-dimensional orbits are identified using a similar argument.  If the 
point $(x_1, x_2,x_3,
y_1, y_2, y_3)$ lies on a one-dimensional orbit then at most two entries
among the $x_i$ and $y_j$ are nonzero.  The nonzero coordinates cannot be
$x_i$ and $y_i$ for the same index $i$ because that point is not on the
hypersurface.  Any other pair of nonzero coordinates defines a one-dimensional orbit
of the torus action.  The fixed points in the closure of each one-dimensional orbit are
the two coordinate points given by the two possible nonzero entries.

We compute the stabilizer of each one-dimensional orbit by examining 
the torus action.
For instance, the 
stabilizer of the orbit for which either $x_1$ and $x_3$ are nonzero 
or $y_1$ and $y_3$ are nonzero is $\C^* \times 1 \times \C^*$.
The other stabilizers can be found by inspection in just the same way.

The graph of fixed points and one-dimensional orbits is given in Figure \ref{moment
graph for hypersurface}.  Each vertex is labelled by the coordinate which
is nonzero, and each edge is labelled by the element of $\t^*$ that
generates
the annihilator of the stabilizer of the orbit.  Note again that
parallel edges have the same label.
\begin{figure}[h]
\begin{picture}(350,120)(0,-60)
\put(175,50){\circle*{7}}
\put(177,-49){\circle*{7}}
\put(160,15){\circle*{7}}
\put(245,15){\circle*{7}}
\put(200,-15){\circle*{7}}
\put(110,-15){\circle*{7}}
\put(174,50){\line(-1,-1){65}}
\put(203,-16){\line(-5,-6){28}}
\put(175,50){\line(2,-1){70}}
\put(110,-15){\line(2,-1){65}}
\put(110,-15){\line(1,0){90}}
\multiput(160,15)(10,0){9}{\line(1,0){5}}
\multiput(107,-16)(15,9){4}{\line(5,3){10}}
\put(200,-16){\line(4,3){44}}
\put(176,50){\line(1,-3){22}}
\multiput(160,14)(4,-16){4}{\line(1,-4){3}}
\put(242,13){\line(-1,-1){63}}
\multiput(174,47)(-8,-20){2}{\line(-2,-5){4}}
\put(170,56){$y_1$}
\put(172,-60){$x_1$}
\put(95, -18){$x_2$}
\put(198, -9){$x_3$}
\put(250,12){$y_2$}
\put(142, 12){$y_3$}
\put(177,-35){$t_2$}
\put(146,-44){$t_3$}
\put(215,-25){$t_1-t_2$}
\put(110,20){$t_1-t_2$}
\put(210,35){$t_3$}
\put(165,30){$t_2$}
\multiput(220,17)(3,1.5){18}{\circle*{1}}
\put(220,17){\line(1,2){5}}
\put(220,17){\line(1,0){11}}
\put(274,41){$t_2-t_3$}
\put(210,0){$t_1$}
\put(133,-1){$t_1$}
\multiput(180,-17)(3,-1.5){20}{\circle*{1}}
\put(180,-17){\line(1,-2){5}}
\put(180,-17){\line(1,0){11}}
\put(240,-48){$t_2 - t_3$}
\multiput(95,-30)(3,0){25}{\circle*{1}}
\put(170,-30){\line(-2,1){10}}
\put(170,-30){\line(-2,-1){10}}
\put(65,-32){$t_1 - t_3$}
\multiput(191,8)(3,0){30}{\circle*{1}}
\put(191,8){\line(2,1){10}}
\put(191,8){\line(2,-1){10}}
\put(282,5){$t_1-t_3$}
\end{picture}
\caption{The graph of fixed points and one-dimensional orbits for a
quadric hypersurface in $\CP^5$} \label{moment graph for hypersurface}
\end{figure}
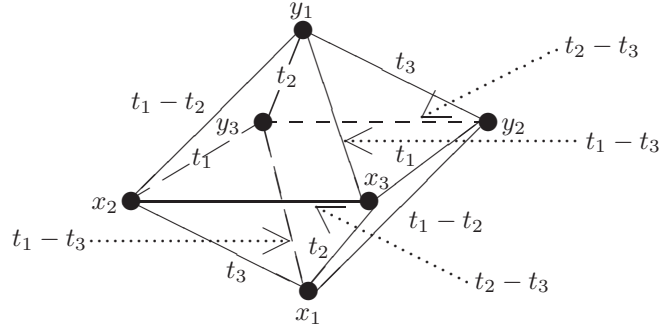

\begin{figure}[h]
\begin{picture}(350,150)(50,-60)
\put(175,50){\circle*{7}}
\put(177,-49){\circle*{7}}
\put(160,15){\circle*{7}}
\put(245,15){\circle*{7}}
\put(200,-15){\circle*{7}}
\put(110,-15){\circle*{7}}
\put(174,50){\line(-1,-1){65}}
\put(203,-16){\line(-5,-6){28}}
\put(175,50){\line(2,-1){70}}
\put(110,-15){\line(2,-1){65}}
\put(110,-15){\line(1,0){90}}
\multiput(160,15)(10,0){9}{\line(1,0){5}}
\multiput(107,-16)(15,9){4}{\line(5,3){10}}
\put(200,-16){\line(4,3){44}}
\put(176,50){\line(1,-3){22}}
\multiput(160,14)(4,-16){4}{\line(1,-4){3}}
\put(242,13){\line(-1,-1){63}}
\multiput(174,47)(-8,-20){2}{\line(-2,-5){4}}
\put(80,72){$\begin{array}{c} 
   p_1+(t_3+t_2-t_1)p_2 +
    t_2(t_2-t_1)p_3 \\ + (t_1-t_3)(t_1-t_2)p_4
   +t_2(t_1-t_2)(t_1-t_3)p_5 \\
   + t_2t_3(t_1-t_2)(t_1-t_3)p_6
  \end{array}$}
\put(172,-60){$p_1$}
\put(63, -18){$p_1+t_3p_2$}
\put(218, -18){$p_1 + t_2p_2 +t_2(t_2-t_3) p_3$}
\put(225,30){$\begin{array}{c}p_1+(t_2-t_1)p_2 +
    (t_2-t_1)(t_2-t_3)p_3 \\ + t_1(t_1-t_2)p_4
   +t_1(t_1-t_2)(t_2-t_3)p_5 \end{array}$}
\put(60, 12){$\begin{array}{l}p_1+(t_3-t_1)p_2 \\+t_1(t_1-t_3)p_4
    \end{array}$}
\end{picture}
\caption{The equivariant cohomology of a quadric hypersurface in $\CP^5$} \label{cohomology of hypersurface}
\end{figure}
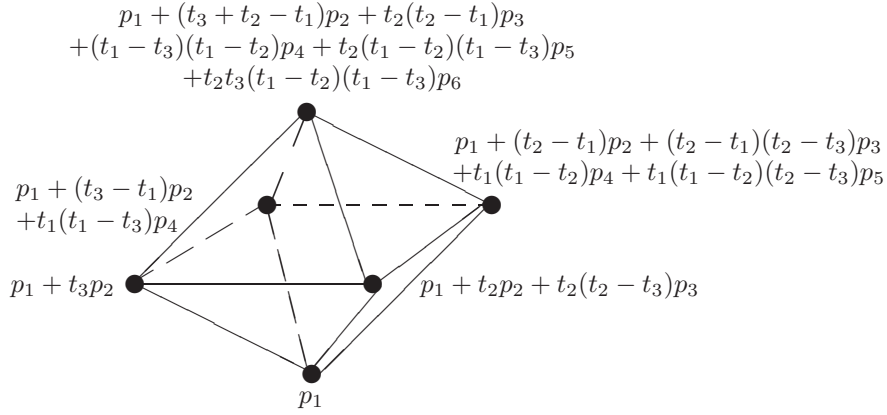

We now use this graph to compute the equivariant cohomology of the 
hypersurface $X$.  Each vertex in Figure \ref{cohomology of hypersurface} 
is labelled by the generic element of
$S(\t^*)$ associated to that fixed point by $H^*_T(X)$.  
It requires no argument to label the first two vertices
$p_1$ and $p_1+t_3p_2$, since no cycles have
yet been formed.  
The vertex for the fixed point $x_3$ can be reached either
by moving over the edge labelled $t_2-t_3$ from the function $p_1 + t_3 p_2$
or up the edge labelled $t_2$ from $p_1$.  Consequently, this vertex is
labelled $p_1 + t_3 p_2 + (t_2-t_3)q$ where $q$ is an element of $S(\t^*)$
that satisfies
$t_3 p_2 + (t_2-t_3)q \in \langle t_2 \rangle$.  Since this means that
$t_3(p_2 - q) \in \langle t_2 \rangle$, we conclude that $q=p_2 + t_2 p_3$
for some $p_3$ freely chosen in $S(\t^*)$.  The general function associated
to the leftmost vertex in the back is identified by a symmetric argument.
The rightmost vertex in the back can be reached by going up the edge
labelled $t_1-t_2$ from $p_1$, going over the edge labelled $t_2-t_3$
from $p_1 + (t_3-t_1)p_2 +t_1(t_1-t_3)p_4$, or going back the edge
labelled $t_1$ from $p_1 + t_2p_2 + t_2(t_2-t_3)p_3$.  This means that the
generic function labelling this vertex is $p_1 + t_2p_2 + t_2(t_2-t_3)p_3
+t_1 q$ for some $q$ in $S(\t^*)$ satisfying both
\[\begin{array}{c}
t_2p_2 + t_2(t_2-t_3)p_3+t_1 q \in \langle t_1 -t_2 \rangle \textup{ and} \\
(t_2p_2 + t_2(t_2-t_3)p_3 +t_1 q) - ((t_3-t_1)p_2 +t_1(t_1-t_3)p_4) 
    \in \langle t_2 - t_3 \rangle. \end{array}\]
To satisfy the first relation, we need $q'$ in $\langle t_1-t_2 
\rangle$ so that $q = - p_2 - (t_2-t_3)p_3 + q'$.  To satisfy the 
second relation, the function $q' = (t_1-t_2)p_4 + (t_1-t_2)(t_2-t_3)p_5$
where $p_5$ is chosen freely in $S(\t^*)$.
The generic function labelling the topmost vertex is shown in the graph.
It is left to the reader to supply the argument, which requires satisfying
the four conditions imposed by the edges leading to the vertex.  
Figure \ref{generator examples} displays one of these generators in a 
very commonly used graphical form.

\begin{example}
The ordinary cohomology ring of Examples \ref{cp2 example}, 
\ref{flag variety example}, and \ref{G24 example}.
\end{example}

Multiplication in the equivariant cohomology ring can be done
coordinatewise.  Since each coordinate of an equivariant cohomology 
class is an element of $S(\t^*)$, this calculation is straightforward.
As an example, we square the class $(0,t_1,t_2)$ in Example 
\ref{cp2 example} and find $(0,t_1,t_2)^2 = (0,t_1^2, t_2^2)$, which
we rewrite in terms of the original classes 
as $t_1 (0,t_1,t_2) + (0,0,t_2(t_2-t_1))$.

Now we use the equality $H^*(X) = \frac{H^*_T(X)}{M \cdot H^*_T(X)}$ from 
Section \ref{defining eq cohom}.  If $u$ represents the equivariant 
cohomology class $(0,t_1,t_2)$ and $v$ represents the class 
$(0,0,t_2(t_2-t_1))$ then our calculation shows that $u^2 = t_1u+v$ in 
equivariant cohomology.  Taking the quotient by $M \cdot H^*_T(\CP^2)$
in effect just sets all $t_i=0$, so in ordinary cohomology $u^2=v$.  This
confirms what we already know to be true.

We leave to the reader the exercise of
computing full multiplication tables
for the ordinary cohomology rings in Examples \ref{cp2 example}, 
\ref{flag variety example}, and \ref{G24 example}.

\section{The moment map and its image} \label{moment map}

The graph of fixed points and one-dimensional orbits of the $T$-action
on $X$ is actually a small part of a larger collection of combinatorial
and geometric information associated to $X$.  We describe this for
a $T$-action on a variety $X$, though it can be defined for general 
$G$-actions by replacing appropriate letters in what follows.
Each complex projective
variety $X$ with an algebraic $T$-action is equipped with a 
moment map $\mu: X \longrightarrow \t^*$, where $\t^*$
denotes the dual to the Lie algebra of $T$, 
or equivalently the cotangent space to $T$
at the identity.  In this section, we describe the moment map in three
different ways: locally, globally, and then combinatorially.  The 
combinatorial description naturally extends the graphs used to
compute the equivariant cohomology of $X$.

We first describe the moment map locally at a point $x$ in $X$.  The 
$T$-action on $X$ induces a linear 
$T$-action on the tangent space $T_xX$,
which we recall 
is diffeomorphic to a neighborhood of $x$ in $X$.  In fact this 
diffeomorphism can be taken to intertwine the local and global 
moment maps.  
Locally near $x$, the 
moment map for $T$ acting on $T_xX$ differs from that for 
$X$ by a translation \cite[Theorem 4.2]{GuSt}.
Since $T_xX$ is a vector space, its moment map is easy to describe
explicitly.  Indeed, $T$ can be viewed as a subgroup of the general 
linear group on $T_xX$.
This permits us to write $\t$ as a subspace of the tangent space to $GL_n$, 
and in particular as some subset of $n \times n$ matrices.  
In this language, the moment map is simply
\[\begin{array}{crcl} \mu: & T_xX & \longrightarrow & \t^* \\
  &y &\mapsto& \{M \mapsto y^TMy\}. \end{array}\]
(This presentation follows that of \cite[Equation 4.5 and preceding]{GuSt}.)

To describe the moment map in general, we make reference to the fact that
a complex projective variety has a naturally associated symplectic form
$\omega$ (see \cite[pages 96--97]{C} for a 
construction of $\omega$, the Fubini-Study form).
Given this form, the moment map is the
$T$-equivariant map $\mu:X \longrightarrow \t^*$ whose
differential satisfies
\begin{equation} \label{moment map defn}
d\mu (\xi)(\alpha) = \omega(\xi, \alpha_x),\end{equation}
where $\xi$ is a tangent vector in $T_xX$, $\alpha$ is in the tangent 
space $\t$ to $T$,
and $\alpha_x$ is the vector field on $X$ induced by the infinitesimal action
along $\alpha$.  (See \cite[Chapter 22]{C} for more.)

For many purposes, all relevant information about the image
of the moment map can be computed without the map itself.  
First, the $T$-fixed points in $X$ will be mapped to points
in $\t^*$.  Let $T'$ be the codimension-one stabilizer of a one-dimensional
orbit $e$ and $\t'$ its tangent.
Applying Equation \ref{moment map defn} to each $\alpha$ in $\t'$ shows 
that the moment map sends the orbit $e$ into a line segment between the 
image of 
two fixed points in $\t^*$.  Moreover, the direction of this line
segment is determined by the annihilator of $\t'$.  
By construction, two line segments are parallel
if their corresponding orbits have the same stabilizer.
Up to rescaling along each edge, then, the closures of the orbits 
determine the image of the moment map.  Lest the reader be too cavalier
about this rescaling, it does correspond to
the choice of a line bundle on the variety.

If $X$ is an irreducible projective variety, 
the image of the moment map is in fact the convex hull
of these fixed points (see \cite{A} and \cite[Theorem 4]{GuSt}).  
In other words, this image 
is a polytope.  
Various classes of varieties have
been studied using the combinatorics of these moment
map images \cite{Gu}, most notably toric varieties \cite{F}.  
From the other direction, 
combinatorial properties of polytopes have been studied by constructing 
projective varieties whose moment map images are the desired polytopes
\cite{St}.  Note that the moment map image of the fixed points 
and one-dimensional orbits
of $X$ are exactly the $0$- and $1$-dimensional faces of this polytope.

Not only can the image of the moment map be reconstructed from the fixed
points
and one-dimensional orbits of $X$ under $T$, but the algebraic 
constructions used in
Section \ref{main theorem} can be extended to any graph satisfying 
certain minimal conditions.  We give the description in \cite{BM} here;
\cite{GuZ} has a similar construction.
A moment graph is defined to be a finite graph for which each edge $e$
is associated to a one-dimensional subspace $V_e$ of the vector space
$\t^*$, and with a partial order on the vertices so that if $u$ and $v$
are joined by an edge then either $u < v$ or $v < u$.  The subspace $V_e$
is called the {\em direction} of the edge $e$.  

In our case, the direction $V_e$ corresponds to the annihilator
of the codimension-one subtorus that stabilizes 
the one-dimensional orbit $e$.  Intuitively, since the closure of
$e$ is 
isomorphic to $\CP^1$, it has both a north pole and a south pole.
We would like to define a partial order so that the north pole is 
greater than the south pole for each one-dimensional orbit.  To do this,
we need to choose a generic map $\rho: \C^* \longrightarrow T$, where 
generic means that the resulting $\C^*$-action on $X$ 
has the same fixed points 
as the original $T$-action.  Given a one-dimensional orbit
$e$ and a point $p$ on $e$, the south pole of $e$ is the fixed point
$v = \lim_{t \to 0} \rho(t) \cdot p$ while the north pole is the fixed
point $u = \lim_{t \to \infty} \rho(t) \cdot p$.
The map $\rho$ in effect chooses a ``global north'' and is discussed
further in the next paragraph.
When $X$ and its $T$-action satisfy certain conditions, for instance if
the $T$-action induces a Whitney stratification \cite{BM} or if the
$T$-action is Hamiltonian \cite[Theorem 1.4.2]{GuZ}, then the
convention that if $u$ is the north pole 
then $u > v$ gives a partial order on the graph.  (Note that there
is an ``opposite'' map that is the composition of $\rho$ 
with the map $z \mapsto \frac{1}{z}$ in $\C^*$.  In terms
of the moment graph, the opposite map 
reverses the original partial order.)

We can also describe this partial order given just the data of
the skeleton of fixed points and one-dimensional orbits
in the moment map image of $X$.
Each edge $e$ in $\t^*$ defines a hyperplane in $\t$,
namely the tangent $\t_e$ of the stabilizer of the corresponding
one-dimensional orbit.  
Since there are a finite number of edges, the complement of the union 
$\bigcup_{e} \t_e$ is dense in $\t$.  In 
particular, we can pick an element $\xi$ in $\t$ so that $e(\xi)$ is
nonzero for each $e$.  
If the endpoints of $e$ are $u$ and $v$, then their difference in $\t^*$
satisfies either $(u-v)(\xi)>0$ or $(v-u)(\xi)>0$.  Define a partial
order on the graph by $u > v$ if $(u-v)(\xi)>0$.  (Transitivity holds since
if $u > v$ and $v > w$ then by linearity $(u-w)(\xi) = (u-v+v-w)(\xi) > 0$.)
This element $\xi$ can also be constructed as 
the derivative of the map $\rho$ mentioned in the previous paragraph.
Since the choice of $\xi$ was only loosely constrained, 
the reader can see that many partial orders can be imposed on the vertices.

In the language of moment graphs, 
Theorem \ref{GKM theorem} associates a copy of $S(\t^*)$ to each vertex 
in the moment graph and identifies $H^*_T(X)$ with the submodule for which
whenever $e$ is an edge with endpoints $u$ and $v$ then 
$f_u - f_v \in S(\t^*)/V_e S(\t^*)$.  This perspective is further studied
in \cite{BM}, which describes more general ways of associating 
$S(\t^*)$-modules to the vertices of moment graphs so that the edges 
impose compatibility restrictions.  In \cite[Section 1.4]{BM}, 
they provide a purely algebraic
algorithm to check these compatibility restrictions and show how this
can be used to construct the equivariant intersection cohomology of $X$.
We note that \cite{BM} assume that $X$ has a
$T$-equivariant stratification by cells, which in turn implies equivariant
formality of the variety.  However, 
their construction depends only on the moment 
graph, regardless of whether it is associated to an algebraic variety or not.
Indeed, some moment graphs 
are the graphs of more than one algebraic variety, not all of which need
be equivariantly formal.  (Compare, for instance, the moment graph of 
the variety shown in Figure \ref{non equivariantly formal} with that of
$\CP^2$.)  Given a fixed moment graph, it is possible
to construct a (noncompact) 
complex manifold with that moment graph so long as some 
added conditions are satisfied.  \cite[Theorem 3.1.1]{GuZ} gives an
explicit construction and a characterization of the extra conditions
needed.  

\section{Generators for each equivariant cohomology group} \label{generators}

In this section we discuss an algorithm to construct a set of generators
for $H^*_T(X)$ as a module over $S(\t^*)=H^*_T(\textup{pt})$.  This method
is implicit in the calculations from the examples of Section \ref{examples}.
We include here both its limitations and some extensions.

Our approach to describing $H^*_T(X)$ in Section \ref{examples}
was to find the most general 
element of $H^*_T(X)$ that could be associated to each fixed point.  
For instance, we began with the lowest vertex and associated a generic 
element in $S(\t^*)$ to it.  We then chose a neighbor $v$ 
and associated the most general element of $S(\t^*)$ that
would satisfy the edge relation to $v$ from the lowest vertex.  We
continued until all vertices were labelled, each time selecting a lowest
neighbor of labelled vertices and determining the most general element of
$S(\t^*)$ that could be associated to that neighbor while satisfying the
edge relations to already labelled vertices.

We now give a computationally simpler and more systematic method for
constructing elements of $H^*_T(X)$ given a moment graph for $X$.
Denote the edge between vertices $w$ and $u$ by $wu$ and denote the
direction of this edge $V_{wu}$.
Fix a vertex $v$ and label all vertices $u$ such that $u < v$ in the 
partial order with the function $0$.  
Now, if $w$ is a vertex for which each
neighbor $u < w$ has been labelled with $f_u \in S(\t^*)$, associate
to $w$ a minimal degree element $f_w \in S(\t^*)$ satisfying
\[f_w - f_u \in V_{wu} \textup{   for all neighbors } u < w.\]
Continue this process until all vertices have been labelled.  This gives
a cohomology class associated to $v$, which we will examine in more detail
later.

It requires some proof even to see that an element of $H^*_T(X)$ can
be produced by the process described.  However, 
this algorithm, which is described in \cite[Section 1.4]{BM}, 
not only works but in fact gives a minimal set of 
generators for $H^*_T(X)$ as a module over $S(\t^*)$ in many cases.  
\cite{GuZ2} gives conditions for the moment graph
which ensure that these cohomology classes actually 
generate the cohomology ring.  

To see how the algorithm works, first 
consider the function associated to $v$ itself 
by this construction.  We indicate the partial order by directing the edges
of the moment graph, so the edges pointing down from $v$ are
labelled $f_1$, $f_2$, $\ldots$, $f_k$ as indicated in Figure \ref{generator
picture}.
\begin{figure}[h]
\begin{picture}(350,50)(0,-30)
\put(175,0){\circle*{5}}
\put(165,-2){$v$}
\put(160,20){\vector(3,-4){14}}
\put(170,20){\vector(1,-4){4}}
\put(175,10){$\cdots$}
\put(205,20){\vector(-3,-2){28}}
\put(175,0){\vector(-3,-2){30}}
\put(175,0){\vector(-1,-4){5}}
\put(175,0){\vector(3,-2){30}}
\put(177,-15){$\cdots$}
\put(138,-15){$f_1$}
\put(160,-20){$f_2$}
\put(208,-15){$f_k$}
\put(143,-30){$0$}
\put(168,-30){$0$}
\put(205,-30){$0$}
\end{picture}
\caption{The edges at the vertex $v$} \label{generator picture}
\end{figure}
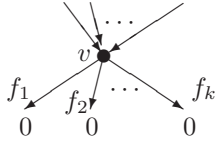
The function associated to $v$ must be
in the ideal generated by the $f_i$.  Since each edge is labelled by
an element of degree one in $S(\t^*)$ and since the lines that they span meet
at a point, the $f_i$ are distinct.  Hence, a natural choice
for the minimum degree symmetric function  
associated to $v$ is $f_1 f_2 \cdots f_k$.  
If the moment graph satisfies the conditions of
\cite[Theorem 2.2]{GuZ2}, then for each fixed point $v$,  
the generator associated to $v$ 
has degree equal to the number of edges pointing downward from $v$.  

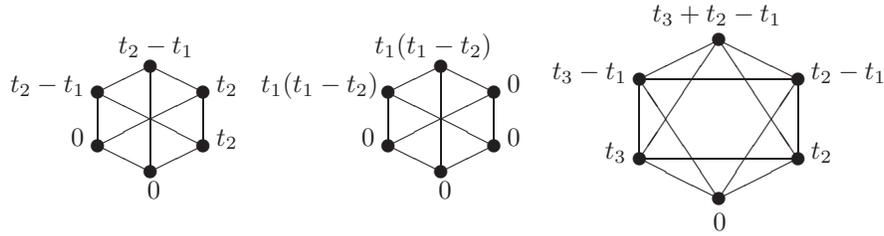
\begin{figure}[h]
\begin{picture}(350,75)(0,-40)
\put(65,20){\circle*{5}}
\put(45,10){\circle*{5}}
\put(85,10){\circle*{5}}
\put(45,-10){\circle*{5}}
\put(85,-10){\circle*{5}}
\put(65,-20){\circle*{5}}
\put(65,20){\line(-2,-1){20}}
\put(65,20){\line(2,-1){20}}
\put(45,-10){\line(0,1){20}}
\put(85,-10){\line(0,1){20}}
\put(65,-20){\line(-2,1){20}}
\put(65,-20){\line(2,1){20}}
\put(45,-10){\line(2,1){40}}
\put(85,-10){\line(-2,1){40}}
\put(65,20){\line(0,-1){40}}
\put(64,-30){$0$}
\put(35,-10){$0$}
\put(90,-10){$t_2$}
\put(12,10){$t_2-t_1$}
\put(90,10){$t_2$}
\put(53, 25){$t_2-t_1$}

\put(175,20){\circle*{5}}
\put(155,10){\circle*{5}}
\put(195,10){\circle*{5}}
\put(155,-10){\circle*{5}}
\put(195,-10){\circle*{5}}
\put(175,-20){\circle*{5}}
\put(175,20){\line(-2,-1){20}}
\put(175,20){\line(2,-1){20}}
\put(155,-10){\line(0,1){20}}
\put(195,-10){\line(0,1){20}}
\put(175,-20){\line(-2,1){20}}
\put(175,-20){\line(2,1){20}}
\put(155,-10){\line(2,1){40}}
\put(195,-10){\line(-2,1){40}}
\put(175,20){\line(0,-1){40}}
\put(174,-30){$0$}
\put(145,-10){$0$}
\put(200,-10){$0$}
\put(107,10){$t_1(t_1-t_2)$}
\put(200,10){$0$}
\put(150, 25){$t_1(t_1-t_2)$}

\put(280,-30){\circle*{5}}
\put(280,30){\circle*{5}}
\put(250,-15){\circle*{5}}
\put(310,-15){\circle*{5}}
\put(250,15){\circle*{5}}
\put(310,15){\circle*{5}}
\put(280,-30){\line(-2,1){30}}
\put(280,-30){\line(2,1){30}}
\put(280,-30){\line(-2,3){30}}
\put(280,-30){\line(2,3){30}}
\put(280,30){\line(-2,-1){30}}
\put(280,30){\line(2,-1){30}}
\put(280,30){\line(-2,-3){30}}
\put(280,30){\line(2,-3){30}}
\put(250,15){\line(1,0){60}}
\put(250,15){\line(0,-1){30}}
\put(310,-15){\line(-1,0){60}}
\put(310,-15){\line(0,1){30}}
\put(278,-42){$0$}
\put(237,-15){$t_3$}
\put(315,-15){$t_2$}
\put(315,15){$t_2-t_1$}
\put(217,15){$t_3-t_1$}
\put(255,37){$t_3+t_2-t_1$}
\end{picture}
\caption{Examples of generators found in Section \ref{examples}}
\label{generator examples}
\end{figure}

Figure \ref{generator examples} 
displays more examples of generators found using this algorithm for
fixed points from examples computed in Section \ref{examples}.  The
generators are the same as those found previously, but the process of
computing them is less laborious.  The way the generators are 
displayed in this figure is very often used by practitioners.  Note that in 
these examples, the generator shown is unique up to rescaling by $\C^*$.

\begin{figure}[h]
\begin{picture}(350,60)(0,-30)
\put(125,20){\circle*{5}}
\put(105,10){\circle*{5}}
\put(145,10){\circle*{5}}
\put(105,-10){\circle*{5}}
\put(145,-10){\circle*{5}}
\put(125,-20){\circle*{5}}
\put(125,20){\line(-2,-1){20}}
\put(125,20){\line(2,-1){20}}
\put(105,-10){\line(0,1){20}}
\put(145,-10){\line(0,1){20}}
\put(125,-20){\line(-2,1){20}}
\put(125,-20){\line(2,1){20}}
\put(123,-32){$0$}
\put(95,-20){$0$}
\put(150,-20){$t_2$}
\put(150,10){$t_2$}
\put(70,10){$t_2-t_1$}
\put(112,25){$t_2-t_1$}

\put(225,20){\circle*{5}}
\put(205,10){\circle*{5}}
\put(245,10){\circle*{5}}
\put(205,-10){\circle*{5}}
\put(245,-10){\circle*{5}}
\put(225,-20){\circle*{5}}
\put(225,20){\line(-2,-1){20}}
\put(225,20){\line(2,-1){20}}
\put(205,-10){\line(0,1){20}}
\put(245,-10){\line(0,1){20}}
\put(225,-20){\line(-2,1){20}}
\put(225,-20){\line(2,1){20}}
\put(223,-32){$0$}
\put(195,-20){$0$}
\put(250,-20){$t_2$}
\put(250,10){$t_1$}
\put(195,10){$0$}
\put(223,25){$0$}
\end{picture}
\caption{A moment graph which is not Palais-Smale, with
two minimal-degree generators for the vertex labelled $t_2$} \label{Palais-Smale counterexample}
\end{figure}
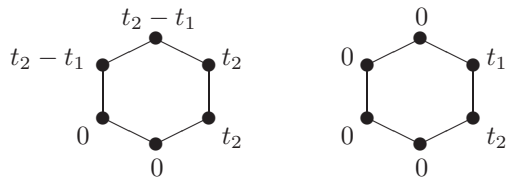

In fact, this algorithm produces a {\em unique} 
set of generators as long as the moment
map for $X$ has a Palais-Smale component.  For our purposes, this 
condition means that the moment graph can be drawn in the plane so that
the vertex $v$ lies above its neighbor $u$ if and only if $v$ has more
downward-pointing edges than $u$.  A formal definition can be found in
\cite{Kn}; see also \cite[Theorem 2.3]{GuZ2}.  
We observe that all of the varieties $X$ discussed
in this paper satisfy this condition.  However, it is easy to find 
varieties which do not.  For instance, recall the flag variety of Example
\ref{flag variety example} and let $S$ be the diagonal matrix with
diagonal entries $1$, $2$, and $3$.  The subvariety of the complex 
flag variety
that consists of flags $V_1 \subsetneq V_2 \subsetneq \C^3$ such that 
$SV_1 \subsetneq V_2$ has the moment graph shown
in Figure \ref{Palais-Smale counterexample}, by a similar argument to
that in Example \ref{flag variety example}.  
By inspection, this moment graph cannot
satisfy the Palais-Smale condition.  We have shown two different 
minimal-degree generators
associated to the same vertex.  (This variety is one of a 
class of subvarieties of the flag variety called Hessenberg varieties;
see \cite{T}.  The cells of the Bia\l ynicki-Birula decomposition, namely the
attracting sets for the torus action, are not in general a 
stratification of Hessenberg varieties.)

A similar approach to that presented here 
can be used to find generators for cohomology rings
related to those of $H^*_T(X)$.  For instance,  in the context of
symplectic manifolds,
 a natural restriction map called the Kirwan map 
sends the equivariant cohomology $H^*_T(X)$ onto the ordinary 
cohomology of the
symplectic reduction of $X$ (which is essentially the
quotient of $X$ by the action of $T$).  \cite{ToW} use this to identify 
a subset of generators for $H^*_T(X)$ that generate the entire 
cohomology $H^*(X /\! / T)$ of the symplectic reduction, and \cite{Go}
constructs these generators combinatorially.

\end{document}